\documentclass{amsart}
\usepackage{amsmath,amssymb,amsfonts,epsfig}

\newtheorem{theorem}{Theorem}[section]

\newtheorem{conjecture}[theorem]{Conjecture}

\theoremstyle{definition}

\theoremstyle{remark}

\numberwithin{equation}{section}

\def\DJ{{\hbox{D\kern-.8em\raise.15ex\hbox{--}\kern.35em}}}
\def\DJo{\DJ okovi\'c}

\def\al{{\alpha}}
\def\be{{\beta}}

\def\de{{\delta}}

\def\vf{{\varphi}}
\def\bC{{\bf C}}

\def\bZ{{\bf Z}}

\def\pI{{\mathcal I}}
\def\pP{{\mathcal P}}

\def\GL{{\mbox{\rm GL}}}
\def\SL{{\mbox{\rm SL}}}
\def\SU{{\mbox{\rm SU}}}

\begin{document}

\title[ Poincar\'{e} series of invariants of binary forms  ]
{A heuristic algorithm for computing the Poincar\'{e} series
of the invariants of binary forms}

\author[D.\v{Z}. \DJo{}]
{Dragomir \v{Z}. \DJo{}}

\address{Department of Pure Mathematics, University of Waterloo,
Waterloo, Ontario, N2L 3G1, Canada}

\email{djokovic@uwaterloo.ca}


\thanks{The author was supported in part by the NSERC
Grant A-5285.}

\date{}

\begin{abstract}
We propose a heuristic algorithm for fast computation of the
Poincar\'{e} series $P_n(t)$ of the invariants of binary forms
of degree $n$, viewed as rational functions. The algorithm is based
on certain polynomial identities which remain to be proved rigorously.
By using it, we have computed the $P_n(t)$ for $n\le30$.
Dixmier proposed long ago three conjectures concerning these
Poincar\'{e} series. We verify that the first two of them
are valid in the above range.
As a supplement to this note, we provide a table from which one
can easily construct the functions $P_n(t)$ for $3\le n\le 30$.
Only the portion of the table covering the range $3\le n\le 20$
is actually displayed in the paper.
\end{abstract}

\maketitle

\section{Introduction}

Let $\pP(V)$ be the algebra of polynomial functions on a
finite-dimensional complex vector space $V$. Let
$\pP^n(V)\subseteq\pP$ be the space of homogeneous polynomials of
degree $n$.
We shall assume from now on that $V=\bC^2$ and refer to
$\pP^n(\bC^2)$ as the space of binary forms of degree $n$.
By fixing a basis, we shall view $\pP(\bC^2)$ as the polynomial
algebra $\bC[x,y]$ in the two coordinate functions $x$ and $y$.
The natural action of $\GL_2(\bC)$ on $V$ induces an action on
$\pP(\bC^2)$ such that each subspace $\pP^n(\bC^2)$ is
stable and the representation of $\GL_2(\bC)$ on it is irreducible.

To simplify the notation, we set $G=\SL_2(\bC)$. One of the main
themes of the classical invariant theory was the study of the algebra of
$G$-invariants of $\pP(\pP^n(\bC^2))$, $n=1,2,3,\ldots$
Denote this subalgebra by $\pI_n=\pP(\pP^n(\bC^2))^G$.
This notation is borrowed from the article of R. Howe \cite{RH},
which gives a modern point of view on this classical subject.
The Poincar\'{e} series (also known as the Hilbert series) of
$\pI_n$ is the formal power series
\[  P_n(t)=\sum_{k\ge0}\dim(\pI_n^k)t^k, \]
where
\[ \pI_n^k = \pI_n\cap\pP^k(\pP^n(\bC^2)). \]

It is well known that this series is in fact the Taylor series
of a rational function, which we denote also by $P_n(t)$.
We consider here the problem of computing efficiently
these rational functions.

There are two methods that we find in the recent literature
on this subject. First there is a method described in
Howe's article and implemented on a computer by P. Saly, Jr.
It is based on the classical formulae, due to Cayley and Sylvester,
for the coefficients of the above Taylor series
and the fact that the concrete form for
the denominator of $P_n(t)$ is known.
Indeed, Dixmier's Conjecture 1 in \cite{JD} gives simple expressions
for the lowest term denominators of the $P_n$'s. (See also
Section \ref{opis} below.)

The second method is based on an explicit but complicated
formula for $P_n(t)$ discovered by T. Springer \cite{TS}.
This formula has been used by Brouwer and Cohen \cite{BC}
and also by Littelmann and Procesi \cite{LP}.

Our method is completely different and is based on certain
conjectural polynomial identities (which are not explicitly
given). In spite of this vagueness, the conjectural formulae
can be used to write an efficient algorithm to compute the
$P_n$'s for small $n$'s. The algorithm saves time by avoiding
the usual tedious procedure of computing many residues
one at the time. We used Maple \cite{Map} to obtain
the formulae for $P_n$ for $n\le30$. The memory requirements
are very modest and the computations can be performed on a PC.

We did not have access to the tables of Brouwer and Cohen
(for $n\le17$) to compare them with our tables.
At the end of his article \cite{RH}, Howe gives
two concrete coefficients of $t^k$
in the numerator of $P_n(t)$, one for $n=19$ and
$k=160$ and the other for $n=20$ and $k=84$.
The denominators are specified as
$\prod_{j=2}^{18}(1-t^{2j})$ for $n=19$, and
$\prod_{j=2}^{19}(1-t^j)$ for $n=20$.
We find that the first coefficient is correct but the second one
should be multiplied by 2.

In a separate file, we provide the list of the $P_n$'s for $n\le30$ as they
may be of interest to other researchers. (I wish I had access to
such a list when I was doing these computations.) More details
about this list are given in Section \ref{opis}.

\section{ The main actors }
\label{uloge}

Let $z$ and $t$ be independent commuting variables and
$\bZ[z,t]$ the corresponding polynomial ring with integer
coefficients.

Let $n\ge3$ be an integer.
Set $n=2s-1$ if $n$ is odd and $n=2s$ if $n$ is even.
The main actors in this paper are the three polynomials
$p_n,q_n\in\bZ[z,t]$ and $r_n\in\bZ[t]$.
The first two are defined, for odd $n$, by
\[ p_n(z,t)=\prod_{i=1}^s (1-tz^{2i-1}),\quad
q_n(z,t)=\prod_{i=1}^s (z^{2i-1}-t), \]
and, for even $n$, by
\[ p_n(z,t)=\prod_{i=1}^s (1-tz^{2i}),\quad
q_n(z,t)=\prod_{i=1}^s (z^{2i}-t). \]

In both cases, $s$ is the $t$-degree of $p_n$ (and $q_n$), i.e.,
its degree as a polynomial in the variable $t$.
Denote by $m$ the $z$-degree of $p_n$ (and $q_n$).
Thus, $m=s^2$ if $n$ is odd and $m=s(s+1)$ if $n$ is even.

Observe that
\begin{equation} \label{p-z-t} 
 z^m p_n(z^{-1},t)=q_n(z,t), \quad (-t)^s p_n(z,t^{-1})=q_n(z,t), 
\end{equation}
and, consequently,
\[  z^m(-t)^s p_n(z^{-1},t^{-1})=p_n(z,t).  \]
These formulae remain valid when the letters $p$ and $q$ are
interchanged.

The third polynomial is defined by
\[ r_n(t)=\prod_{i=2}^{n-1} (1-t^{2i}) \]
if $n$ is odd, and by
\[ r_n(t)=(1+t)\prod_{i=2}^{n-1} (1-t^i) \]
if $n$ is even.

Finally, set
\[ \vf_n(z,t)=z^{m-2} (z^2-1) r_n(t). \]
It is easy to verify that
\begin{equation} \label{f-z-t} 
 z^{2m-2}\vf_n(z^{-1},t)=-\vf_n(z,t), \quad
 t^d \vf_n(z,t^{-1})=(-1)^n \vf_n(z,t), 
\end{equation}
and, consequently,
\[  z^{2m-2}t^d \vf_n(z^{-1},t^{-1})=(-1)^{n+1}\vf_n(z,t),  \]
where $d$ is the degree of $r_n(t)$. Thus $d=2s(n-2)$ for $n$ odd
and $d=s(n-1)$ for $n$ even.

\section{ Basic conjecture }

In this section we state our basic working conjecture, which will be
taken for granted in the remaining part of the paper.
It is supported by extensive computer calculations.

Let $I_n=<p_n,q_n>$ be the ideal of $\bZ[z,t]$ generated by
$p_n$ and $q_n$.
We conjecture that $\vf_n\in I_n$. More precisely,

\begin{conjecture} \label{hip-1}
There exist unique polynomials $a_n,b_n\in\bZ[z,t]$ of $z$-degree
$m-2$ such that
\begin{equation} \label{f-a-b}
\vf_n(z,t)=a_n(z,t)q_n(z,t)+b_n(z,t)p_n(z,t).
\end{equation}
\end{conjecture}

(The uniqueness is clear because $p_n$ and $q_n$ are relatively
prime.)

Let us give an example. If $n=5$ then $s=3$, $m=9$ and $\vf = aq+bp$,
where
\begin{eqnarray*}
\vf &=& z^7(z^2-1)(1-t^4)(1-t^6)(1-t^8), \\
p &=& (1-tz)(1-tz^3)(1-tz^5), \\
q &=& (z-t)(z^3-t)(z^5-t), \\
a &=& (1-t^6+t^{12})-t(1+t^2-t^6)(1-t^2-t^6)z+t^2(1-t^6+t^8)z^2 \\
&& -t(1-t^4)(1-t^6-t^8)z^3+t^2(1-t^4)(1+t^2-t^8)z^4-t^5(1-t^2+t^8)z^5 \\
&& +t^2(1+t^4-t^6)(1-t^4-t^6)z^6-t^3(1-t^6+t^{12})z^7, \\
b &=& t^3(1-t^6+t^{12})-t^2(1+t^4-t^6)(1-t^4-t^6)z+t^5(1-t^2+t^8)z^2 \\
&& -t^2(1-t^4)(1+t^2-t^8)z^3+t(1-t^4)(1-t^6-t^8)z^4-t^2(1-t^6+t^8)z^5 \\
&& +t(1+t^2-t^6)(1-t^2-t^6)z^6-(1-t^6+t^{12})z^7.
\end{eqnarray*}

By applying the substitution $z\to z^{-1}$ to (\ref{f-a-b}) and by using the
identities (\ref{p-z-t}) and (\ref{f-z-t}), we obtain that
\begin{equation} \label{for-b1}
b_n(z,t)=-z^{m-2} a_n(z^{-1},t).
\end{equation}
Similarly, the substitution $t\to t^{-1}$ gives
\begin{equation} \label{for-b2}
b_n(z,t)=(-1)^{n+s} t^{d-s} a_n(z,t^{-1}).
\end{equation}
From the last two equations we deduce that
\begin{equation} \label{for-a}
z^{m-2}t^{d-s} a_n(z^{-1},t^{-1})=(-1)^{n+s+1} a_n(z,t).
\end{equation}

We can write
\begin{equation} \label{alfa-beta}
a_n(z,t)=\sum_{i=0}^{m-2} \al_k(t)z^k, \quad
b_n(z,t)=\sum_{i=0}^{m-2} \be_k(t)z^k,
\end{equation}
where $\al_k,\be_k\in\bZ[t]$.
Then the equation (\ref{for-b1}) gives
\begin{equation} \label{for-beta}
\be_k(t)=-\al_{m-2-k}(t), \quad 0\le k\le m-2
\end{equation}
and (\ref{for-a}) gives
\begin{equation} \label{for-alfa}
\al_{m-2-k}(t)= (-1)^{n+s+1} t^{d-s}\al_k(t^{-1}), \quad 0\le k\le m-2.
\end{equation}

Hence, the equation (\ref{f-a-b}) can now be rewritten as
\begin{equation} \label{glavna}
\vf_n(z,t)=\sum_{k=0}^{m-2} \al_k(t)
\left( z^k q_n(z,t)-z^{m-2-k} p_n(z,t) \right).
\end{equation}

By setting $z=0$ in (\ref{glavna}), we obtain that
$\al_{m-2}(t)=(-t)^s \al_0(t)$. By combining this with
the equation (\ref{for-alfa}) for $k=0$, we obtain that
\begin{equation} \label{for-0}
\al_0(t)=(-1)^{n+1} t^{d-2s} \al_0(t^{-1}).
\end{equation}

Apparently the computation of the polynomials
$a_n(z,t)$ is  a difficult job. Fortunately, we do not
need to know these polynomials explicitly but only the
polynomials $a_n(0,t)=\al_0(t)$. The main point of our
algorithm is that the polynomials $\al_0(t)$ can be
computed efficiently.

For the curious reader, let us throw in one more intriguing
conjecture (which plays no role whatsoever in this paper).

\begin{conjecture} \label{hip-2}
$I_n\cap\bZ[t]$ is the principal ideal of $\bZ[t]$ generated by
the polynomial
\[ (1-t^2)\prod_{i=1}^{n-1} (1-t^{2i}),\quad
(1+t)\prod_{i=1}^{n-1} (1-t^i),\quad
\prod_{i=1}^{n-1} (1-t^i), \]
according to whether $n$ is odd, congruent to $2$ modulo $4$,
or divisible by $4$.
\end{conjecture}
This conjecture makes sense (and is true) also when $n$ is
equal to 1 or 2.

\section{The integral formula}

The Molien--Weyl integral formula for $P_n(t)$ can be written
in the following form (see \cite{DK})
\[ P_n(t)=\frac{1}{2\pi i} \int_{|z|=1} f_n(z,t)
\frac{ {\rm d}z }{z}, \]
where
\[ f_n(z,t)=\frac{1-z^{-2}}{\prod_{k=0}^n (1-tz^{n-2k})}. \]
The integration is to be performed over the unit circle in the
counterclockwise direction and it is assumed that $|t|<1$.

Let us define yet another function
\[ g_n(z,t)=\frac{z^{m-2}(z^2-1)}{p_n(z,t)q_n(z,t)}. \]
For odd $n$, $f_n(z,t)=g_n(z,t)$ and, for even $n$,
$f_n(z,t)=g_n(z,t)(1-t)^{-1}$. Since the
integration variable is $z$, the factor $(1-t)^{-1}$ can be
inserted at the very end of the computation.
Thus it suffices to compute the integral
\[ Q_n(t)=\frac{1}{2\pi i} \int_{|z|=1} g_n(z,t)
\frac{ {\rm d}z }{z}. \]

The formula (\ref{f-a-b}) gives
\[ r_n(t)g_n(z,t)=\frac{a_n(z,t)}{p_n(z,t)}+
\frac{b_n(z,t)}{q_n(z,t)}. \]
Hence,
\[ r_n(t)Q_n(t)=
\frac{1}{2\pi i} \int_{|z|=1} \frac{a_n(z,t)}{p_n(z,t)}
\frac{ {\rm d}z }{z}+
\frac{1}{2\pi i} \int_{|z|=1} \frac{b_n(z,t)}{q_n(z,t)}
\frac{ {\rm d}z }{z}. \]
The second integral is 0 since all the poles are inside the
unit circle (and the residue at $\infty$ is 0).
The integrand of the first integral has all of its poles
outside the unit circle except for the simple pole at $z=0$.
Since $a_n(0,t)=\al_0(t)$ and $p_n(0,t)=1$, we obtain that
\begin{equation} \label{for-Q}
Q_n(t)=\frac{\al_0(t)}{r_n(t)}.
\end{equation}

Hence, the computation of $P_n(t)$ is reduced to computing
the polynomial $\al_0(t)$.
Since $Q_n(0)=r_n(0)=1$, we deduce that $\al_0(0)=1$.
The formula (\ref{for-0}) now implies that $\al_0$ has
degree $d-2s$. This is in agreement with the well known
fact (it was known to Sylvester \cite{SF}) that the rational function
$P_n(t)$ has degree $-n-1$.

If $n$ is odd, then $\al_0(t)$ is a palindromic polynomial, i.e.,
its coefficients are symmetric about the midpoint.
If $n$ is even, then $\al_0(1)=0$, i.e., $\al_0(t)$ is divisible
by $1-t$.

Let us emphasize that the main feature of our algorithm described above
is that we compute the above integral at one fell swoop, avoiding
the tedious procedure of computing all the residues for the poles
inside the unit circle one at the time.

\section{How to compute $\al_0$ ?}

Let me start by quoting the Rule \#1 from the list of useful tips
on Jean-Charles Faug\`{e}re's home page: ``In a first time try
to compute modulo $p$ where $p$ is a small prime.''
I learned this rule on my own (i.e., the hard way).

To compute $\al_0(t)$ we used the Maple package called 
``LinearAlgebra'' and its subpackage called ``Modular''.
First we choose a big prime, $l$, in the data type
integer[4]. The largest one is $l=65521$. We perform computations
modulo this prime. We generate random integer mod $l$ inputs
for the variable $z$ and plug it into our equation (\ref{glavna}).
We need only $m-1$ inputs since there are $m-1$
unknown $\al_k$'s.

Thus we obtain a system of $m-1$ equations which are linear
in the $\al_k$'s. But these are still polynomial equations
as they contain the variable $t$. We chose a random value,
say $\tau$ modulo $l$,
for the variable $t$ and plug it into this system
of equations. Now we have just a small system of linear
equations to solve for the unknowns $\al_k(\tau)$. This job
can be easily handled by ``Modular''. We are only interested
in the value $\al_0(\tau)$.

Next recall that we know the degree of $\al_0$; it is equal
to $d-2s$. Thus we have to repeat the above calculation
$d-2s+1$ times to get the values $\al_0(\tau)$ for
$d-2s+1$ different inputs $\tau$. Having done this,
the $d-2s+1$ unknown coefficients of $\al_0(t)$ can be easily
computed by solving the corresponding Van der Monde system of
linear equations.
Of course, this means that at this stage we know these
coefficients modulo our prime $l$.

Next we have to repeat the whole calculation above with
several distinct primes. We used the largest seven primes
available: 
\[ 65521,\, 65519,\, 65497,\, 65479,\, 65449,\, 65447,\, 65437. \]
This is necessary when $n$ is near the high end of our range because
some of the coefficients of $\al_0(t)$ for $n=29$ have 30 digits.
The remaining task is to lift these modular solutions to
a genuine solution over $\bZ$. For this purpose, ``Modular''
provides the ``ChineseRemainder'' command which makes
the task very easy.

In the hardest case, $n=29$, the computation of $\al_0(t)$ modulo
the above seven primes (performed one after the other) took
almost 5 hours on the SunBlade workstation running a single
R10000 CPU at 250 MHz with 8 GB of RAM.

Certainly, if Conjecture \ref{hip-1}  was false these computations,
which involve so many random choices, would not produce a
meaningful result. Moreover, the polynomials $\al_0(t)$ have
certain symmetry properties.
As mentioned earlier, if $n$ is odd then $\al_0(t)$ is a
palindromic polynomial. If $n$ is even, then $\al_0(t)$ is divisible
by $1-t$, according to one of the Dixmier's conjectures.
Finally, if $n$ is divisible by 4, then $\al_0(t)$ is
also divisible by $1+t^{s-1}$, again by Dixmier.
This is all true in the cases that we have computed.

If $n\ge3$ then $t=1$ is a pole of $P_n(t)$ of order $n-2$ and
the limit of $(1-t)^{n-2}P_n(t)$ as $t\to1$ has been computed by
Hilbert \cite{DH}. We have verified that our formulae for $P_n(t)$
agree with his result.

\section{Description of the coefficient table} \label{opis}

There exist unique relatively prime polynomials $A_n,B_n\in\bZ[t]$,
with $A_n(0)=B_n(0)=1$, such that
\[ P_n(t)=\frac{A_n(t)}{B_n(t)}. \]
In this notation, Dixmier's Conjecture 1 asserts that
$B_n(t)=r_n(t)/(1+t^{s-1})$ if $n$ is divisible by 4, and
$B_n(t)=r_n(t)$ otherwise. He verified this conjecture for $n\le16$
and our computations show that it is valid for $n\le30$.

We propose the following additional conjecture:
\begin{conjecture} \label{hip-3}
The polynomials $A_n(t)$ are irreducible for $n\ge5$.
\end{conjecture}

Since the denominators $B_n(t)$ are easy to write down, we
provide only the information necessary to construct the
numerators $A_n(t)$ for $3\le n\le30$.

If $n$ is odd, then $A_n(t)=\al_0(t)$.
Now assume that $n$ is even. Then 
$\al_0(t)=(1-t)\bar{\al}_0(t)$ for some $\bar{\al}_0\in\bZ[t]$ and
$P_n(t)=\bar{\al}_0(t)/r_n(t)$. Thus
$A_n(t)=\bar{\al}_0(t)$ if $n\equiv2 \pmod{4}$.
Finally, if $n\equiv0 \pmod{4}$, then
$\bar{\al}_0(t)$ is divisible by $1+t^{s-1}$,
and we have
$A_n(t)=\bar{\al}_0(t)/(1+t^{s-1})$.

The degree of $A_n$ is always even and we shall denote it by $2\de$.
We have $\de=2s(s-2)$ for $n$ odd, $2\de=s(2s-3)-1$ for
$n\equiv2 \pmod{4}$, and $\de=s(s-2)$ for $n\equiv0 \pmod{4}$.
These numerators are palindromic polynomials
\[ A_n(t)=c_0+c_1t+\cdots+c_{\de-1} t^{\de-1}+c_\de t^\de
+c_{\de-1} t^{\de+1}+\cdots+c_1 t^{2\de-1}+c_0 t^{2\de} \]
and so it suffices to store roughly half of their coefficients.

For each $n$ in the range $3\le n\le 30$ we record the coefficients
$c_0=1,c_1,\ldots,c_\de$ of $A_n$ in that order.
They are stored as a Maple table $A[3..30]$.

For example, if $n=5$ then $s=3$, $\de=6$, $A_5=1-t^6+t^{12}$.
This is recorded in our table as 
\[ A[5]=[1,0,0,0,0,0,-1]. \]

Assume that $n\ge3$ is odd. Then it is known that $c_i=0$ for odd
$i$'s as well as for $i=2$. Dixmier \cite{JD} shows that
$c_4$ and $c_8$ are nonnegative, and so are $c_6$ for
$n\ge15$ and $c_{10}$ for $n\ge9$. We can prove that in fact
\[ c_4=\left[ \frac{n-3}{6} \right]. \]
It may be useful to obtain similar formulae for $c_6$, $c_8$, etc.

\section{Appendix: Table $A[3..20]$} \label{tablica}

We display here the coefficients $c_0=1,\; c_1,\ldots,\; c_\delta$ of
the numerator $A_n(t)$. Recall that the degree of $A_n(t)$ is $2\delta$.
The display includes only the values $n=3,4,\ldots,20$. The full table
with the range $n=3,4,\ldots,30$ is available from the author on request.

\newpage
\begin{eqnarray*}
 && A[3] = [1] \\
 && A[4] = [1] \\
 && A[5] = [1,0,0,0,0,0,-1] \\
 && A[6] = [1,1,0,-1,-1] \\
 && A[7] = [1,0,0,0,0,0,-1,0,2,0,-1,0,5,0,2,0,6] \\
 && A[8] = [1,1,0,-1,-1,0,1,1,1] \\
 && A[9] = [1,0,0,0,1,0,-1,0,5,0,3,0,18,0,15,0,44,0,43,0,82,0,76,0,122,0, \\
&& 107,0,147,0,119] \\
 && A[10] = [1,1,0,-1,-1,-1,2,2,4,6,10,10,16,15,16,16,19,15] \\
 && A[11] = [1,0,0,0,1,0,-1,0,10,0,11,0,58,0,85,0,222,0,336,0,660,0,951,0, \\
&& 1589,0,2154,0,3188,0,4080,0,5510,0,6633,0,8310,0,9443,0,11059,0, \\
&& 11894,0,13094,0,13319,0,13852] \\
 && A[12] = [1,1,0,0,1,1,3,7,11,15,25,35,50,70,95,121,156,191,229,266,304, \\
&& 332,358,371,377] \\
 && A[13] = [1,0,0,0,1,0,-1,0,19,0,31,0,157,0,321,0,885,0,1756,0,3794,0, \\
&& 6856,0,12788,0,21324,0,35633,0,55326,0,85174,0,124064,0,178645,0, \\
&& 246238,0,334814,0,439321,0,568305,0,712862,0,881834,0,1061455,0, \\
&& 1259989,0,1459221,0,1666984,0,1860904,0,2049854,0,2209072,0,2349306, \\
&& 0,2446352,0,2514111,0,2530530] \\
 && A[14] = [1,1,0,-1,0,0,5,8,21,34,69,107,193,295,477,705,1064,1490,2121, \\
&& 2865,3876,5040,6535,8186,10217,12379,14898,17490,20381,23197,26216, \\
&& 28992,31799,34192,36461,38127,39559,40263,40636] \\
 && A[15] = [1,0,0,0,2,0,0,0,32,0,76,0,378,0,995,0,3048,0,7294,0,17681,0, \\
&& 37736,0,78903,0,152321,0,285968,0,507762,0,876759,0,1451423,0,2341739, \\
&& 0,3653241,0,5568497,0,8254649,0,11983447,0,16987847,0,23631274,0, \\
 && 32196429,0,43116834,0,56681420,0,73342055,0,93320393,0,117007543,0, \\
 && 144461993,0,175919353,0,211175615,0,250222591,0,292516508,0, \\
&& 337751801,0,385016863,0,433713649,0,482605505,0,530877973,0, \\
&& 577086324,0,620343376,0,659172312,0,692798202,0,719914717,0, \\
&& 740045690,0,752239053,0,756462172] \\
 && A[16] = [1,1,0,0,1,3,10,18,38,79,153,278,514,891,1523,2528,4072,6367, \\
&& 9772,14572,21306,30498,42785,58863,79666,105852,138459,178275,226114, \\
&& 282612,348514,423785,508764,603037,705993,816528,933532,1054746, \\
&& 1178503,1302160,1423129,1538529,1645908,1741955,1824693,1891546, \\
&& 1940709,1970683,1980976]
\end{eqnarray*}

\begin{eqnarray*}
 && A[17] = [1,0,0,0,2,0,0,0,50,0,156,0,844,0,2716,0,9280,0,26055,0, \\
&& 70846,0,173224,0,405183,0,883551,0,1847356,0,3669433,0,7024773,0, \\
&& 12919848,0,23019526,0,39697193,0,66608244,0,108748704,0,173371011, \\
&& 0,270001994,0,411791616,0,615371715,0,902700319,0,1300556398,0, \\
&& 1842885348,0,2569619659,0,3529481145,0,4777707107,0,6379225544,0, \\
&& 8404807944,0,10934524315,0,14051849433,0,17847385164,0,22410522390, \\
&& 0,27833575972,0,34200807232,0,41593316021,0,50075498973,0, \\
&& 59701312756,0,70498510402,0,82477120228,0,95612069584,0, \\
&& 109854817368,0,125115124389,0,141276752124,0,158179073390,0, \\
&& 175637520454,0,193425539961,0,211299997259,0,228983624510,0, \\
&& 246195050087,0,262631400199,0,278002047047,0,292010632283,0, \\
&& 304391087240,0,314888097678,0,323293686822,0,329425050810,0, \\
&& 333160013567,0,334411422423] \\
 && A[18] = [1,1,0,-1,1,2,11,20,60,122,292,573,1199,2264,4307,7692,13639, \\
 && 23121,38688,62619,99808,154969,236962,354532,522744,756614,1080091, \\
 && 1517149,2103843,2875718,3883971,5178127,6826649,8894132,11466985, \\
 && 14623600,18466440,23083747,28588450,35071007,42645944,51394256, \\
&& 61419462,72778864,85549508,99748823,115411657,132499497,150990759, \\
&& 170778600,191774520,213796608,236686776,260188811,284080692, \\
&& 308043009,331804545,355002792,377341521,398442572,418014893, \\
&& 435697784,451237581,464324429,474773647,482354991,486977592, \\
&& 488512665] 
\end{eqnarray*}

\newpage
\begin{eqnarray*}
 && A[19] = [1,0,0,0,2,0,1,0,76,0,296,0,1763,0,6738,0,25712,0,82893,0, \\
&& 252012,0,694765,0,1807368,0,4392969,0,10154779,0,22296093,0, \\
&& 46930799,0,94802787,0,184822672,0,348299741,0,636837951,0,1131559371, \\
&& 0,1959027689,0,3309329549,0,5465457626,0,8835569146,0,14002953513,0, \\
&& 21778902561,0,33281420196,0,50015102190,0,73986433683,0,107815040622, \\
&& 0,154890725357,0,219515780752,0,307103168571,0,424341546490,0, \\
&& 579425240340,0,782222455286,0,1044516156680,0,1380143077777,0, \\
&& 1805211986577,0,2338168202449,0,2999945060393,0,3813898895206,0, \\
&& 4805843429850,0,6003797685274,0,7437855144286,0,9139694108669,0, \\
&& 11142247265979,0,13478942963820,0,16183166152183,0,19287230767266, \\
&& 0,22821658964018,0,26813925195712,0,31287624164380,0,36261071290807, \\
&& 0,41746472475910,0,47748507080547,0,54263646677747,0,61278874674037, \\
&& 0,68771319597259,0,76707272862777,0,85042295486194,0,93720678352722, \\
&& 0,102676146969267,0,111831862933442,0,121101792448393,0, \\
&& 130391278393710,0,139599058259333,0,148618354053166,0, \\
&& 157339436459345,0,165651094224795,0,173443559767368,0, \\
&& 180610146679903,0,187050288671933,0,192671089114197,0, \\
&& 197390193869982,0,201136986494180,0,203855025614692, \\
&& 0,205502649265697,0,206054755643582]
\end{eqnarray*}

\begin{eqnarray*}
 && A[20] = [1,1,0,0,2,5,17,40,100,232,544,1199,2599,5365,10770,20867, \\
&& 39312,71826,128004,222286,377375,626606,1019690,1627231,2550571, \\
&& 3929105,5955664,8888409,13072800,18958777,27131530,38333433,53503793, \\
&& 73805811,100673389,135839979,181391167,239789025,313924970,407134237, \\
&& 523239787,666546890,841871660,1054502755,1310204347,1615137665, \\
&& 1975827494,2399034646,2891684840,3460685330,4112816074,4854497808, \\
&& 5691649536,6629418196,7672023583,8822464074,10082371655,11451725459, \\
&& 12928740940,14509620528,16188514360,17957332700,19805806095, \\
&& 21721392176,23689449970,25693256426,27714311797,29732461463, \\
&& 31726325603,33673521588,35551188212,37336279422,39006154571, \\
&& 40538894181,41913908026,43112232559,44117095575,44914144829, \\
&& 45491928155,45842000355,45959277535]
\end{eqnarray*}


\begin{thebibliography}{99}

\bibitem{BC}
A.E. Brouwer and A.M. Cohen, \emph{The Poincar\'{e} series of the 
polynomial invariants under $\SU_2$ in its irreducible representation
of degree $\le17$}, preprint of the Mathematisch Centrum,
Amsterdam, 1979.

\bibitem{DK}
H. Derksen and G. Kemper, \emph{Computational Invariant Theory},
Springer-Verlag, New York, 2002.

\bibitem{JD}
J. Dixmier, \emph{Quelques r\'{e}sultats et conjectures concernant
les s\'{e}ries de Poincar\'{e} des invariants des formes binaires},
in ``Lecture Notes in Mathematics,'' Vol. 1146, pp. 127--160,
Springer-Verlag, Berlin, 1985.

\bibitem{DH}
D. Hilbert, \emph{\"{U}ber die vollen Invariantensysteme}, Math. Ann.
\textbf{36} (1890), 473--534.

\bibitem{RH}
R. Howe, \emph{The Classical Groups} and Invariants of Binary Forms,
in \emph{The Mathematical Heritage of Hermann Weyl},
Vol. 48, Proc. Symp. Pure Math., Amer. Math. Soc., Providence,
R.I., 1987.

\bibitem{LP}
P. Littelmann and C. Procesi, \emph{On the Poincar\'{e} series
of the invariants of binary forms}, J. Algebra \textbf{133} (1990),
490--499.

\bibitem{Map}
MAPLE, {\em Maplesoft}, Waterloo, Ontario.

\bibitem{TS}
T.A. Springer, \emph{On the invariant theory of $\SU_2$},
Indag. Math. {\bf 42} (1980), 339-345.

\bibitem{SF}
J.J. Sylvester and F. Franklin, \emph{Tables of the generating functions
and groundforms for the binary quantics of the first ten orders},
Amer. J. Math. {\bf 2} (1879), 223-251.

\end{thebibliography}
\end{document}